\documentclass[11pt]{article}
\usepackage[short]{optidef}
\usepackage{tikz, lscape, amsmath,amsfonts,latexsym,amssymb, amsthm}
\usepackage[utf8]{inputenc}
\usepackage[english]{babel}
\usepackage{lmodern, verbatim}
\usepackage[T1]{fontenc}
\usepackage{calc,bm}
\usepackage{tabularx,enumitem,nameref,hyperref}
\usepackage[mathscr]{eucal}
\usepackage{subfig,stackrel,floatrow}

\renewcommand{\b}{\mathbf}

\usepackage[numbers]{natbib}
\usepackage[margin=1.5cm]{geometry}%[bottom=80pt]

\newcommand{\NN}{\mathbb{N}}

\newcommand{\RR}{\mathbb{R}}

\newcommand{\Hk}{\mathcal{H}_K}

\newcommand{\Sx}{\mathcal{S}_{[t_0,T]}}
\newcommand{\Kx}{K_{[t_0,T]}}

\newcommand{\T}{\mathscr{T}}

\DeclareMathOperator*{\Img}{Im}

\newcommand{\Id}{\text{Id}}
\renewcommand{\d}{\mathrm{d}} %dt
\newcommand{\tb}{\textbf} %bold text
\newcommand{\R}{\RR} %reals
\renewcommand{\b}{\mathbf}
\newcommand{\rv}{}%\textcolor{blue}
 
\newcommand{\iv}[2]{[\![#1,#2]\!]} %interval: [[a,b]]

\theoremstyle{definition}

\theoremstyle{plain}
\newtheorem{Theorem}{Theorem}

\newtheorem{Lemma}{Lemma}%[Theorem]

\theoremstyle{definition}
\newtheorem{Definition}{Definition}%[Theorem]
\date{\today}

\makeatletter
\let\orgdescriptionlabel\descriptionlabel
\renewcommand*{\descriptionlabel}[1]{%
	\let\orglabel\label
	\let\label\@gobble
	\phantomsection
	\edef\@currentlabel{#1\unskip}%
	\let\label\orglabel
	\orgdescriptionlabel{#1}%
}
\makeatother

\title{Interpreting the dual Riccati equation through the LQ reproducing kernel}
\author{Pierre-Cyril Aubin-Frankowski\footnote{\'Ecole des Ponts ParisTech and CAS, MINES ParisTech, PSL Research University, France. Email: pierre-cyril.aubin@mines-paristech.fr}}
\date{\today}%

\begin{document}
	\maketitle
	\begin{abstract}
		In this study, we provide an interpretation of the dual differential Riccati equation of Linear-Quadratic (LQ) optimal control problems. Adopting a novel viewpoint, we show that LQ optimal control can be seen as a regression problem over the space of controlled trajectories, and that the latter has a very natural structure as a reproducing kernel Hilbert space (RKHS). The dual Riccati equation then describes the evolution of the values of the LQ reproducing kernel when the initial time changes. This unveils new connections between control theory and kernel methods, a field widely used in machine learning.
	\end{abstract}
	\section{Introduction}
	We consider the problem of finite-dimensional time-varying linear quadratic (LQ) optimal control with finite horizon and quadratic terminal cost as in
	\begin{mini!}
		{\substack{\b u(\cdot)}}{\b x(T)^\top \b J_T \, \b x(T)+\int_{t_0}^{T}[\b x(t)^{\top}\b Q(t) \b x(t) + \b u(t)^{\top} \b R(t) \b u(t)] \d t}{\label{opt-LQR}}{V(t_0,\b x_0):= \tag{1}}
		\addConstraint{ \b x'(t)}{=\b A(t) \b x(t) + \b B(t) \b u(t), \, \text{a.e.\ in} \, [t_0,T] }{\label{eq_dynamic}}
		\addConstraint{ \b x(t_0)}{=\b x_0,}{ \label{eq_CI}}
	\end{mini!}
	where the state $\b x(t) \in \R^N$ and the control $\b u(t) \in \R^M$. We shall henceforth assume that $\b J_T\succ \b 0$,\footnote{Here $\succcurlyeq$ (resp. $\succ$) denotes the (strict) partial order over positive semi-definite matrices.} and for all $t\in[t_0,T]$, $\b R(t) \succcurlyeq r \Id_M$ with $r>0$, as well as $\b A(\cdot)\in L^1([t_0,T],\R^{N,N})$, $\b B(\cdot)\in L^2([t_0,T],\R^{N,M})$, $\b Q(\cdot)\in L^1([t_0,T],\R^{N,N})$, and $\b R(\cdot)\in L^2([t_0,T],\R^{N,N})$. To have a finite objective, we restrict our attention to measurable controls satisfying $\b R(\cdot)^{1/2} \b u(\cdot)\in  L^2([t_0,T],\R^{N})$. Problem \eqref{opt-LQR} is intimately related to the differential Riccati equation,\footnote{The index $T$ in $\b J(\cdot,T)$ is kept as a reminder that \eqref{eq_Riccati} is defined w.r.t.\ a given terminal time $T$. We denote by $\partial_1$ the derivative w.r.t.\ the first variable.} expressed as
	\begin{equation} \label{eq_Riccati}
	- \partial_1\b J(t,T)= \b A(t)^\top \b J(t,T) + \b J(t,T)\b A(t) - \b J(t,T) \b B(t) \b R(t)^{-1} \b B(t)^\top \b J(t,T) + \b Q(t)\; ; \; \b J(T,T) = \b J_T,
	\end{equation}%\citep[e.g.][pp.\ 31, 408]{Bensoussan2007}
	which solution $\b J(\cdot,T)$ satisfies $V(t_0,\b x_0)=\b x_0^\top \b J(t_0,T) \b x_0$. It is well-known (e.g.\ \cite[pp.\ 31, 408]{Bensoussan2007}) that under the above positivity assumptions, $\b J(t,T)$ is a symmetric positive definite matrix, which inverse $\b M(t,T):=\b J(t,T)^{-1}$ satisfies a dual Riccati equation
	\begin{equation} \label{eq_dual_Riccati}
	\partial_1\b M(t,T)= \b A(t) \b M(t,T) + \b M(t,T)\b A(t)^\top - \b B(t) \b R(t)^{-1} \b B(t)^\top + \b M(t,T)\b Q(t) \b M(t,T)\; ; \; \b M(T,T) = \b J_T^{-1}.
	\end{equation}
	This inverse matrix has been used as a tool to obtain a representation formula in infinite-dimensional LQ control \citep{barbu1992} but it has not received the deserved interest yet. Whereas the solution of \eqref{eq_Riccati} is equal to the Hessian of the value function $V(t_0,\cdot)$, i.e.\ $\b J(t_0,T)=\partial_{\b x,\b x}V(t_0,\cdot)$, we show (Theorem \ref{thm_Riccati-LQkernel} below) that the solution of \eqref{eq_dual_Riccati} is equal to the diagonal element of a matrix-valued reproducing kernel $K(\cdot,\cdot)$, naturally associated with \eqref{opt-LQR}. Owing to this interpretation, the dual Riccati equation \eqref{eq_dual_Riccati} is thus no less fundamental and effectively allows to reverse the perspective between the adjoint vector and the optimal trajectory.\\
	
%	\newP{Luenberger?}
	%	In \citep{aubin2020hard_control}, we considered the inner product \eqref{def_K-normbis} rather than \eqref{def_K-norm}. This was a more appropriate choice for fixed initial time and general lower semicontinuous convex terminal costs, alongside affine state constraints. 
	We first need to bring trajectories to the fore in \eqref{opt-LQR}. In his seminal book, Luenberger \citep[p255]{luenberger1968optimization} already discussed that an optimal control problem such as \eqref{opt-LQR} can be seen as either optimizing over the set of controls $\b u(\cdot)$, or jointly over the set of trajectories $\b x(\cdot)$ and controls $\b u(\cdot)$, connected through the dynamic constraint \eqref{eq_dynamic}. Luenberger also alluded without details to a third possibility, that of optimizing directly over the controlled trajectories. We follow this last viewpoint and consequently introduce the vector space $\Sx$ of controlled trajectories  of the linear system:
	\begin{equation}
	\Sx:=\{\b x:[t_0,T]\rightarrow\R^N \,|\, \exists \, \b u(\cdot) \text{ s.t.\ } \b x'(t)=\b A(t) \b x(t) + \b B(t) \b u(t) \text{ a.e.\ and $\b R(\cdot)^{1/2} \b u(\cdot)\in  L^2([t_0,T],\R^{N})$} \}.\label{def_Sx}
	\end{equation} 
	There is not necessarily a unique choice of $\b u(\cdot)$ for a given $\b x(\cdot) \in \Sx$.\footnote{This is the case for instance if $\b B(t)$ is not injective for a set of times $t$ with positive measure.} Therefore, with each $\b x(\cdot) \in \Sx$, we associate the control $\b u(\cdot)$ having minimal norm based on the pseudoinverse $\b B(t)^{\ominus}$ of $\b B(t)$ for the $\R^M$-norm $\|\cdot\|_{\b R(t)}:=\|\b R(t)^{1/2}\cdot\|$:
	\begin{align}
	\b u(t) = \b B(t)^{\ominus}[\b x'(t) - \b A(t) \b x(t)] \, \text{ a.e.\ in} \,[t_0,T]. \label{def_u_as_X-X'} %L(t)([\b x(t),\b x'(t)]) 
	\end{align}
	Problem \eqref{opt-LQR} then induces a natural inner product over $\Sx$. As a matter of fact, the expression
	\begin{equation}
	\left<\b x_1(\cdot),\b x_2(\cdot)\right>_{K}:= \b x_1(T)^\top \b J_T\, \b x_2(T) + \int_{t_0}^{T}[\b x_1(t)^\top \b Q(t)\b x_2(t) + \b u_1(t)^\top \b R(t)\b u_2(t)] \d t 
	\label{def_K-norm}
	\end{equation}
	is bilinear and symmetric over $\Sx \times \Sx$. It is positive definite over $\Sx$ as $\|\b x(\cdot)\|^2_K = \b 0$ implies that $\b u(\cdot)\stackrel{a.e}{\equiv} 0$ and, as $\b J_T\succ \b 0$, $\b x(T)= \b 0$, hence $\b x(\cdot)\equiv \b 0$. Therefore
	\begin{mini}
		{\substack{\b x(\cdot)\in\Sx}}{\|\b x(\cdot)\|^2_K}{\label{opt-LQ_Gram}}{V(t_0,\b x_0)=}%\tag{$\Psc_K$}
		\addConstraint{ \b x(t_0)}{=\b x_0}.
	\end{mini}
	In other words the value function $V(t_0,\b x_0)$ of \eqref{opt-LQR} coincides with the optimal value of a constrained norm minimization over $\Sx$. The solution of \eqref{opt-LQ_Gram} can be made explicit as $(\Sx,\left<\cdot,\cdot\right>_{K})$ is not an arbitrary Hilbert space, but a vector-valued reproducing kernel Hilbert space (vRKHS).
	
	\section{Vector spaces of linear controlled trajectories as vRKHSs}	\label{sec_RKHS}
	\begin{Definition}\label{def_vRKHS}
		Let $\T$ be a non-empty set. A Hilbert space $(\Hk(\T),\left<\cdot,\cdot\right>_{K})$ of $\R^N$-vector-valued functions defined on $\T$ is called a vRKHS if there exists a matrix-valued kernel $K_\T:\T \times \T \rightarrow \R^{N,N}$ such that the \emph{reproducing property} holds: for all $t \in \T,\, \b p\in\R^N $,  $K_\T(\cdot,t)\b p \in \Hk(\T)$ and for all $\b f \in \Hk(\T)$, $\b p^{\top}\b f(t) = \left<\b f,K_\T(\cdot,t)\b p\right>_{K}$.
	\end{Definition}
\noindent \tb{Remark:} It is well-known that by Riesz's theorem, an equivalent definition of a vRKHS is that, for every $t \in \T$ and $\b p\in\R^N$, the evaluation functional $\b f\in\Hk(\T) \mapsto \b p^\top \b f(t)\in\RR$ is continuous. There is also a one-to-one correspondence between the kernel $K_\T$ and the vRKHS $(\Hk(\T),\left<\cdot,\cdot\right>_{K})$ (see e.g. \citep[Theorem 2.6]{micheli_matrix-valued_2014}). Moreover, by symmetry of the scalar product, the matrix-valued kernel has a Hermitian symmetry, i.e.\ $K_\T(s,t)=K_\T(t,s)^\top$ for any $s,t\in\T$.
	\begin{Lemma}\label{lem_SxHilbert} $(\Sx,\left<\cdot,\cdot\right>_{K})$ is a vRKHS over $[t_0,T]$ with a reproducing kernel $\Kx$ which we call the \emph{LQ kernel}.
	\end{Lemma}
	\noindent \tb{Proof of Lemma \ref{lem_SxHilbert}}: The proof is identical to the one of Lemma 1 in \citep{aubin2020hard_control} where $\Sx$ was equipped with the following inner product
	 \begin{equation}
	\left<\b x_1(\cdot),\b x_2(\cdot)\right>_{K,init}:= \b x_1(t_0)^\top \b x_2(t_0) + \int_{t_0}^{T}[\b x_1(t)^\top \b Q(t)\b x_2(t) + \b u_1(t)^\top \b R(t)\b u_2(t)] \d t.
	\label{def_K-normbis}
	\end{equation}
	\begin{flushright}
		$\blacksquare$
	\end{flushright}

	% By Riesz's theorem, an equivalent definition of a vRKHS is that, for every $t \in \T$ and $\b p\in\R^N$, the evaluation functional $\b f\in\Hk \mapsto \b p^\top \b f(t)\in\RR$ is continuous. Moreover, for any vRKHS $(\Hk,\left<\cdot,\cdot\right>_{K})$, its reproducing kernel $K(\cdot,\cdot)$ is unique. We will identify later the reproducing kernel of $(\Sx,\left<\cdot,\cdot\right>_{K})$. By definition, vRKHSs handle easily pointwise evaluations. Besides, 
	Owing to Lemma \ref{lem_SxHilbert}, we can look for a ``representer theorem``, i.e.\ a necessary condition to ensure that the solutions of an optimization problem like \eqref{opt-LQ_Gram} live in a finite dimensional subspace of $\Sx$ and consequently enjoy a finite representation.
	\begin{Theorem}\citep{aubin2020hard_control} \label{thm_representer} Let $(\Hk(\T),\left<\cdot,\cdot\right>_{K})$ be a vRKHS defined on a set $\T$. For a given $I\in\NN$, let $\{t_{i}\}_{i\in \iv{1}{I}}\subset \T$. Consider the following optimization problem with ``loss`` function $L:\RR^{I}\rightarrow \RR\cup\{+\infty\}$, strictly increasing ``regularizer`` function $\Omega:\RR_+\rightarrow \RR$, and vectors $\{\b c_{i}\}_{i\in \iv{1}{I}} \subset\R^N$
		\begin{mini*}
			{\substack{\b f\in\Hk(\T)}}{ L\left(\b c_{1}^\top \b f(t_{1}),\dots, \b c_{I}^\top \b f(t_{I})\right) + \Omega\left(\|\b f \|_K\right).}{}{}
		\end{mini*}
	Then, for any minimizer $\bar{\b f}$, there exists $\{\b p_{i}\}_{i\in \iv{1}{I}}\subset \RR^N$ such that $\bar{\b f}= \sum_{i=1}^{I} K_\T(\cdot,t_{i}) \b p_{i}$ with $\b p_{i}=\alpha_{i} \b c_{i}$ for some $\alpha_{i}\in\RR$.
	\end{Theorem}
	Taking $L(\b e_1^\top \b x(t_0),\dots,\b e_N^\top \b x(t_0)):=\chi_{\b x_0}(\b x(t_0))$ and $\Omega(y)=y^2$, with $\b e_i$ the $i$-th basis vector of $\R^N$, $\chi_{\b x_0}$ the indicator function of $\b x_0$, we apply Theorem \ref{thm_representer} to \eqref{opt-LQ_Gram}. Since $\|\cdot\|^2_K$ is strongly convex and there exists $\b x(\cdot)\in\Sx$ satisfying $\b x(t_0)=\b x_0$, the solution of \eqref{opt-LQ_Gram} is unique and can be written as $\bar{\b x}(t)=\Kx(t,t_0) \b p_0$, with $ \b p_0=\Kx(t_0,t_0)^\ominus \b x_0 \in\R^N$, where $\Kx(t_0,t_0)^\ominus$ is the pseudoinverse of $\Kx(t_0,t_0)$ for the $\R^N$-seminorm $\|\Kx(t_0,t_0)^{1/2}\cdot\|$. Thus, owing to the reproducing property,
	\begin{align} \label{eq_value_function_kernel}
	V(t_0,\b x_0)=\|\bar{\b x}(\cdot)\|^2_K &= \langle \Kx(\cdot,t_0) \b p_0, \Kx(\cdot,t_0) \b p_0 \rangle_K = \b p_0^\top \Kx(t_0,t_0) \b p_0\\
	&= \b x_0^\top \Kx(t_0,t_0)^{\ominus} \b x_0 = \b p_0^\top \b x_0. \nonumber
	\end{align}
	 So we conjecture that $\Kx(t_0,t_0)^{\ominus}=\b J(t_0,T)$. We actually have a stronger result:
	
	\begin{Theorem}\label{thm_Riccati-LQkernel} Let $K_d:t_0 \in]-\infty,T]\mapsto \Kx(t_0,t_0)$. Then $K_d(t_0)=\b J(t_0,T)^{-1}$.
		%Let $K_d:t_0\mapsto \Kx(t_0,t_0)$ be the map of the diagonal initial values of the reproducing kernel associated with $(\Sx,\left<\cdot,\cdot\right>_{K})$ when varying $t_0\in]-\infty,T]$. Then $\b J(t_0,T)^{-1}=K_d(t_0)$ where $\b J(\cdot,T)$ is the solution of the differential Riccati equation \eqref{eq_Riccati}.
	\end{Theorem}%\b M(t_0,T)
	The proof of Theorem \ref{thm_Riccati-LQkernel} (in Section \ref{sec_proof_Thm2} below) boils down to identifying the reproducing kernel of $(\Sx,\left<\cdot,\cdot\right>_{K})$. Informally, the inverse relation comes from inverting the graph of the $(\b x, \b p)$-relation. As a matter of fact, consider the solution $\b p(t)$ of the adjoint equation
	\begin{equation} \label{eq_adjoint_vect}
	\b p'(t)= - \b A(t)^\top \b p(t) + \b Q(t) \bar{\b x}(t) \quad \quad \b p(T) = -\b J_T\, \bar{\b x}(T).
	\end{equation}
	Then we have $\b p(t)= -\b J(t,T) \bar{\b x}(t)$. In other words, the solution $\b J(\cdot,T)$ of the differential Riccati equation maps the optimal trajectory $\bar{\b x}(\cdot)$ to its adjoint vector $\b p(\cdot)$. On the contrary, since $\bar{\b x}(t)=\Kx(t,t_0) \b p_0$, the kernel $\Kx(\cdot,t_0)$ maps an initial covector $\b p_0 \in \R^n $ to the optimal trajectory $\bar{\b x}(\cdot)$. This effectively inverts the graph of the relation between $\bar{\b x}(\cdot)$ and $\b p(\cdot)$. The inversion performed is related to yet another change of perspective, from an online and differential approach to an offline and integral one.\\ 
	
	Through Pontryagine's Maximum Principle (PMP), it is well known that the optimal control $\bar{\b u}(\cdot)$ satisfies $\bar{\b u}(t)=\b R(t)^{-1} \b B(t)^\top \b p(t)=-\b R(t)^{-1} \b B(t)^\top \b J(t,T) \bar{\b x}(t)=: \b G(t) \bar{\b x}(t)$. Hence, based on $\b J(t,T)$, one has a closed feedback loop, with gain matrix $\b G(t)$, and knows the control to apply based only on the present time and state. However the optimal trajectory $\bar{\b x}(\cdot)$ is not encoded as simply as in the kernel formula $\bar{\b x}(t)=\Kx(t,t_0) \b p_0$. It has to be derived through numerical approximations of the dynamics \eqref{eq_dynamic}. Conversely, the kernel $\Kx$ performs the integration of the Hamiltonian system \eqref{eq_dynamic}-\eqref{eq_adjoint_vect} and sparsely encodes $\bar{\b x}(\cdot)$ over $[t_0,T]$ by $\b p_0$. This sparsity partly stems from the smaller number of constraints in \eqref{opt-LQ_Gram} w.r.t.\ \eqref{opt-LQR} since the dynamics \eqref{eq_dynamic} were incorporated in the definition of $\Sx$. Unlike in the PMP, the adjoint vector $\b p(t)$ disappears in the kernel perspective and only the initial condition (or some intermediate rendezvous points) induce a covector $\b p_i$.
	
	More generally, for a given interval $[t_0,T]$, Theorem \ref{thm_representer} states that to encode the optimal trajectories one needs at most as many covectors $\b p_i$ as there are points $t_i$ where the trajectory is evaluated in the optimization problem. It is a classical property of ``kernel machines``, frequently leveraged in classification tasks (e.g.\ SVMs in \cite{scholkopf2002learning}). This result was exploited in \citep{aubin2020hard_control} to tackle affine state constraints. From the PMP perspective, it resulted in focusing only on the measures supported on the constraint boundary. Unlike the adjoint vector $\b p(t)$ associated with the equality constraint \eqref{eq_dynamic}, which never vanishes except for abnormal trajectories, the covectors corresponding to inequality constraints are null whenever the constraint is not active. This led to extremely sparse encoding of the optimal trajectory by specifying only the active covectors on the $[t_0,T]$ time interval.\footnote{	In \citep{aubin2020hard_control}, we considered the inner product \eqref{def_K-normbis} rather than \eqref{def_K-norm} which assumes a terminal quadratic cost. The choice of \eqref{def_K-normbis} was more appropriate for fixed initial time and general lower semicontinuous convex terminal costs, alongside affine state constraints.} Offline computation of the kernel is indeed well suited for path-planning problems. The kernel formalism however conflicts with the online perspective since varying $t_0$ changes the domain of $\Kx$. As the correspondence between the kernel $K_\T$ and the vRKHS $(\Hk(\T),\left<\cdot,\cdot\right>_{K})$ is one-to-one (e.g.\ \citep[Theorem 2.6]{micheli_matrix-valued_2014}), varying $\T=[t_0,T]$ or modifying the inner product changes the kernel. In general, restricting the domain leads to complicated relations between a vRKHS and its kernel \cite[pp.78-80]{saitoh16theory}. In our case, the dual Riccati equation \eqref{eq_dual_Riccati} precisely describes how the values of the LQ kernel change when varying $t_0$.

	\section{Proof of Theorem \ref{thm_Riccati-LQkernel}}\label{sec_proof_Thm2}
	The proof corresponds to the identification of the reproducing kernel of $(\Sx,\left<\cdot,\cdot\right>_{K})$. Since we shall proceed with fixed initial time $t_0$, we drop the corresponding index and set $K(\cdot,\cdot)=\Kx(\cdot,\cdot)$. By existence and unicity of the reproducing kernel, we just have to exhibit a function $K(\cdot,\cdot)$ which satisfies the requirements of Definition \ref{def_vRKHS}.
	
	Let us denote by $\b \Phi_{\b A}(t,s)\in\R^{N,N}$ the state-transition matrix of $\b z'(\tau)=\b A(\tau) \b z(\tau)$, defined from $s$ to $t$, i.e.\ $\b z(t)=\b \Phi_{\b A}(t,s)\b z(s)$. The key property used throughout this section is \rv{the variation of constants, a.k.a.\ Duhamel's principle,} stating that for any absolutely continuous $\b x(\cdot)$ such that $ \b x'(t)=\b A(t) \b x(t) + \b B(t) \b u(t)$ a.e., we have for any $\sigma,t \in [t_0,T]$
	\begin{align}
	\b x(t) &= \b \Phi_{\b A}(t,\sigma) \b x(\sigma)+\int_{\sigma}^{t} \b \Phi_{\b A}(t,\tau) \b B(\tau) \b u(\tau) \d \tau.\label{eq_Duhamel}
	\end{align}
	Setting $\partial_1 K(s,t):\b p\mapsto \frac{d}{ds} (K(s,t)\b p)$, let us define formally $\b U(s,t):= \b B(s)^{\ominus}[\b \partial_1 K(s,t) - \b A(s) K(s,t)]$. The reproducing property for $K$ then writes as follows, for all $t \in [t_0,T]$, $\b p\in\R^N$, $\b x(\cdot)\in\Sx$,
	\begin{equation}
	\b p^{\top}\b x(t)=(K(T,t)\b p)^{\top} \b J_T \b x(T) + \int_{t_0}^{T} (K(s,t)\b p)^{\top} \b Q(s) \b x(s) \d s+\int_{t_0}^{T} (\b U(s,t)\b p))^{\top} \b R(s) \b u(s) \d s. \label{eq_repro_Q}
	\end{equation}
	By the Hermitian symmetry of $K$ and \rv{the variation of constants} \eqref{eq_Duhamel} written for $\sigma=T$, we can rewrite \eqref{eq_repro_Q} as, \rv{for all $t \in [t_0,T]$, $\b x(\cdot)\in\Sx$,}
	\begin{align*}
	\b x(t)&=K(t,T) \b J_T \b x(T) +  \int_{t_0}^{T} K(t,s) \b Q(s) \left(\b \Phi_{\b A}(s,T) \b x(T)+\int_{T}^{s} \b \Phi_{\b A}(s,\tau) \b B(\tau) \b u(\tau)  \d \tau\right) \d s\\
	&\hspace{3cm}+\int_{t_0}^{T}  \b U(s,t)^{\top} \b R(s) \b u(s) \d s.
	\intertext{Regrouping terms,}
	 \b x(t)&=\left(K(t,T) \b J_T + \int_{t_0}^{T} K(t,s) \b Q(s) \b \Phi_{\b A}(s,T) \d s\right)\b x(T) + \int_{t_0}^{T} K(t,s) \b Q(s) \int_{T}^{s} \b \Phi_{\b A}(s,\tau) \b B(\tau) \b u(\tau)  \d \tau \d s\\
	 &\hspace{3cm}+\int_{t_0}^{T}  \b U(s,t)^{\top} \b R(s) \b u(s) \d s
	 \intertext{Setting $ \tilde{K}(t,s):=\int_{t_0}^{s} K(t,\tau)  \b Q(\tau) \b \Phi_{\b A}(\tau,s)  \d \tau$, and applying Fubini's theorem, we get}
	 \b x(t)&=\left(K(t,T) \b J_T + \tilde{K}(t,T) \right)\b x(T) + \int_{t_0}^{T} \left[\b U(s,t)^{\top} \b R(s)-\int_{t_0}^{s} K(t,\tau) \b Q(\tau)  \b \Phi_{\b A}(\tau,s)\b B(s)\d \tau\right]  \b u(s)   \d s.
	\end{align*}
	Identifying with \eqref{eq_Duhamel} for $\sigma=T$, we derive that
	\begin{align}
	K(t,T)\b J_T=& \; \b \Phi_{\b A}(t,T)-\tilde{K}(t,T),\nonumber \\
	\b U(s,t)^{\top} \b R(s)=&	\left\{\begin{array}{cl}
	(-\b \Phi_{\b A}(t,s)+\tilde{K}(t,s))\b B(s)\; &\forall s\geq t,\\ 
	\tilde{K}(t,s)\b B(s)\; &\forall s < t.
	\end{array} \right.\label{eq_integralK}
	\end{align}
	Let us introduce formally an adjoint equation defined for a variable $\b \Pi(s,t)\in\RR^{N,N}$. For any given $t\in[t_0,T]$,
	\begin{align*}
		\partial_1 \b \Pi(s,t)=- \b A(s)^\top \b \Pi(s,t) +\b  Q(s)K(s,t) \quad \quad \b  \Pi(t_0,t)=-\Id_N.
	\end{align*}
	This is the matrix version of \eqref{eq_adjoint_vect} but with an initial rather than terminal condition. Again, applying \rv{the variation of constants} \eqref{eq_Duhamel} with $\sigma=t_0$ to $\b \Pi(s,t)$, taking the transpose and owing to the symmetries of $K(\cdot,\cdot)$ and $\b \Phi(\cdot,\cdot)$, we derive that
	\begin{align*}
	\b \Pi(s,t) &= \b \Phi_{(-\b A^\top)}(s,t_0) \b \Pi(t_0,t)+\int_{t_0}^{s} \b \Phi_{(-\b A^\top)}(s,\tau) \b  Q(\tau)K(\tau,t)  \d \tau\\
	\b \Pi(s,t)^\top &= -\b \Phi_{\b A}(t_0,s)+\int_{t_0}^{s}  K(t,\tau) \b  Q(\tau) \b \Phi_{\b A}(\tau,s)  \d \tau= -\b \Phi_{\b A}(t_0,s)+\tilde{K}(t,s).
	\end{align*}
	Since $\partial_1 K(s,t) = \b A(s) K(s,t) + \b B(s) \b U(s,t)$, by \eqref{eq_integralK}, for any given time $t\in[t_0,T]$, we have two coupled differential equations for $K(\cdot,t)$ and $\b \Pi(\cdot,t)$ 
	\begin{equation}\label{eq_Hamilton_full}
	\left.\begin{array}{ll}
	\partial_1 K(s,t) &= \b A(s) K(s,t) + \b B(s) \b R(s)^{-1} \b B(s)^\top \left\{
	\begin{array}{cl}
	\b \Pi(s,t) + \b \Phi_{\b A}(t_0,s)^\top - \b \Phi_{\b A}(t,s)^\top\;\forall s\geq t,\\ 
	\b \Pi(s,t) + \b \Phi_{\b A}(t_0,s)^\top\quad\forall s < t.
	\end{array} \right.\\ 
	\partial_1 \b \Pi(s,t) &=- \b A(s)^\top \b \Pi(s,t) +\b  Q(s)K(s,t),\\
	\b  \Pi(t_0,t)&=-\Id_N \quad \quad K(t,T)\b J_T=-\b  \Pi(T,t)^\top +  \b \Phi_{\b A}(t,T) -  \b \Phi_{\b A}(t_0,T).  	\end{array}\right.
	\end{equation}
	Equations \eqref{eq_Hamilton_full} seem quite intricate but they become simpler for $t=t_0$, and as seen in \eqref{eq_value_function_kernel}, $t_0$ is actually the only time that interests us to solve the LQ optimal control problem \eqref{opt-LQ_Gram}. For $t=t_0$, the equations boil down to
	\begin{equation}\label{eq_Hamilton_short}
	\left.\begin{array}{ll}
	\partial_1 K(s,t_0) &= \b A(s) K(s,t_0) + \b B(s) \b R(s)^{-1} \b B(s)^\top \b \Pi(s,t_0)\\ 
	\partial_1 \b \Pi(s,t_0) &=- \b A(s)^\top \b \Pi(s,t_0) +\b  Q(s)K(s,t_0) \\
	\b  \Pi(t_0,t_0)&=-\Id_N \quad ; \quad \b  \Pi(T,t_0)=-\b J_T K(T,t_0). \end{array} \right.
	\end{equation}
	Let us solve \eqref{eq_Hamilton_short} by variation of the constant $\b J_T$. We thus look for a function  $\b J(\cdot,T)$, which we will prove solves \eqref{eq_Riccati}, such that $\b J(T,T)=\b J_T$ and $ \b \Pi(s,t_0)= - \b J(s,T) K(s,t_0)$. We take the derivative in $s$ of the latter expression to obtain
	\begin{align*}
	- \b A(s)^\top \b \Pi(s,t_0) +\b  Q(s)K(s,t_0)  &=- \b J(s,T) \left(\b A(s) K(s,t_0) + \b B(s) \b R(s)^{-1} \b B(s)^\top \b \Pi(s,t_0)\right) - (\partial_1 \b J(s,T))K(s,t_0)%\\&\hspace{2cm}
%	\\ \b A(s)^\top \b J(s,T) K(s,t_0)+ \b  Q(s)K(s,t_0) &= -\b J(s,T) \b A(s) K(s,t_0) + \b J(s,T) \b B(s) \b R(s)^{-1} \b B(s)^\top \b J(s,T) K(s,t_0)\\
%	&\hspace{2cm} - (\partial_1 \b J(s,T))K(s,t_0)
	\end{align*}
	Therefore, applying to the right of the equations any pseudo-inverse of $K(s,t_0)$,
	\begin{align*}
	\b 0=\left[\b A(s)^\top \b J(s,T) ^\top+ \b J(s,T) \b A(s) -\b J(s,T) \b B(s) \b R(s)^{-1} \b B(s)^\top \b J(s,T)  +\b  Q(s) + \partial_1 \b J(s,T)\right] \text{proj}_{\Img K(s,t_0)}.
	\end{align*}
	So it suffices that $\b J(\cdot,T)$ solves the differential Riccati equation \eqref{eq_Riccati} and that, by symmetry of $K$,  
	$$ \b \Pi(t_0,t_0)=-\Id_N= - \b J(t_0,T) K(t_0,t_0) = - K(t_0,t_0) \b J(t_0,T).$$
	Consequently $K(t_0,t_0)=\b J(t_0,T)^{-1}$ which formalizes our intuition of the inverse relation between $\b J$ and $K$. Now let us vary the initial time and consider the function $K_d:t_0\mapsto \Kx(t_0,t_0)$. Taking the derivative w.r.t.\ $t_0$ of $ \b \Pi(t_0,t_0)=-\Id_N= - \b J(t_0,T) \Kx(t_0,t_0)$, we get
	\begin{align*}
	\b 0&=(\partial_1 \b J(t_0,T)K_d(t_0) +\b J(t_0,T)(\partial_1K_d(t_0)),
	\intertext{applying $K_d(t_0)$ to the left of the equation and using that $K_d(t_0)=\b J(t_0,T)^{-1}$, we obtain}
	\b 0&=-K_d(t_0)\b A(t_0)^\top-\b A(t_0)K_d(t_0) +\b B(t_0) \b R(t_0)^{-1} \b B(t_0)^\top - K_d(t_0) \b  Q(t_0) K_d(t_0) + \partial_1K_d(t_0).
	\end{align*}
	This concludes our proof as $K_d(\cdot)$ solves the dual matrix Riccati equation \eqref{eq_dual_Riccati} which has a unique solution.
	
%	\setcitestyle{numbers}
	\bibliographystyle{unsrtnat}
	\bibliography{LQR_Riccati.bib} 

\end{document}